\documentclass[11pt]{article}
\pdfoutput=1

\usepackage[utf8]{inputenc}
\usepackage[T1]{fontenc}
\usepackage[backend=biber,
            bibstyle=phys,
            bibencoding=utf8,
            texencoding=utf8,
            sorting=none,
            url=false,
            doi=false,
            eprint=true,
            sortcites=true,
            citestyle=numeric-comp]{biblatex}
\DeclareFieldFormat*{title}{\mkbibitalic{#1}}
\DeclareFieldFormat{labelnumberwidth}{\mkbibbrackets{#1}~~}
\renewbibmacro*{volume+number+eid}{%
  \printfield{volume}%
  \iffieldundef{number}{}{%
  \printtext{(}%
  }%
  \printfield{number}%
  \iffieldundef{number}{}{%
  \printtext{)}%
  }%
  \setunit{\addcomma\space}%
  \printfield{eid}}
\usepackage{amsfonts,amsthm,amsmath,amssymb}
\usepackage[cal=txupr,
            calscaled=1.,
            scr=esstix,
            scrscaled=1.1,
            bb=txof]{mathalpha} 
\usepackage{latexsym,xcolor}
\usepackage{pifont}
\usepackage{tikz}
\usepackage[all]{xy}
\usepackage[breaklinks]{hyperref} 
\usepackage{authblk}

\theoremstyle{plain}
\newtheorem{theo}{Theorem}[section] 
\newtheorem{prop}[theo]{Proposition}

\newtheoremstyle{rem}{3pt}{3pt}{}{}{\bfseries}{.}{.5em}{}
\theoremstyle{rem}
\newtheorem{remax}[theo]{Remark}
\newenvironment{rema}
{\pushQED{\qed}\remax}
{\popQED\endremax}
\newtheorem{examx}[theo]{Example}
\newenvironment{exam}
{\pushQED{\qed}\examx}
{\popQED\endexamx}
\theoremstyle{definition}
\newtheorem{defi}[theo]{Definition}

\def\sC{{\mathscr{C}}}
\def\sD{{\mathscr{D}}}

\def\cA{{\mathcal{A}}}
\def\cB{{\mathcal{B}}}
\def\cV{{\mathcal{V}}}
\def\cW{{\mathcal{W}}}

\def\bn{{\mathbf{n}}}
\def\be{{\mathbf{e}}}
\def\bx{{\mathbf{x}}}
\def\by{{\mathbf{y}}}
\def\bi{{\mathbf{i}}}
\def\bj{{\mathbf{j}}}
\newcommand{\Bell}[0]{\ensuremath{\boldsymbol\ell}}

\newcommand{\pFq}[5]{{
{}_{#1}F_{#2}\left( \genfrac{}{}{0pt}{0}{#3}{#4}\middle| #5\right)
}}
\newcommand{\rac}[3]{{
R\left( \genfrac{}{}{0pt}{0}{#1}{#2}\middle| #3\right)
}}
\newcommand{\hah}[3]{{
H\left( \genfrac{}{}{0pt}{0}{#1}{#2}\middle| #3\right)
}}

\newcommand{\qqquad}[0]{{\quad\qquad}}
\newcommand{\sss}{{\scriptscriptstyle}}

\usepackage{geometry}
\title{\bf Change of basis for the tridiagonal pairs of type II}
\renewcommand*{\Affilfont}{\normalsize\small}
\author[1a]{Nicolas Cramp\'e\,}
\author[2b]{Julien Gaboriaud\,}
\author[2c]{Satoshi Tsujimoto\,\vspace{.5em}}

\affil[1]{Institut Denis-Poisson CNRS/UMR 7013 - Universit\'e de Tours -
Universit\'e d'Orl\'eans,
\newline\vspace{.9em}
Parc de Grandmont, 37200 Tours, France.}

\affil[1]{Laboratoire d'Annecy-le-Vieux de Physique Théorique LAPTh,
Universit\'e Savoie Mont Blanc,
\newline\vspace{.9em}
CNRS, F-74000 Annecy, France.}

\affil[2]{\vspace{.9em}Graduate School of Informatics, Kyoto University,
Sakyo-ku, Kyoto, 606-8501, Japan.
}
{
 \makeatletter
 \renewcommand\AB@affilsepx{: \protect\Affilfont}
 \makeatother
 \affil[ ]{E-mail addresses}
 \makeatletter
 \renewcommand\AB@affilsepx{, \protect\Affilfont}
 \makeatother
 \affil[a]{crampe1977@gmail.com}
 \affil[b]{julien.gaboriaud@umontreal.ca}
 \affil[c]{tsujimoto.satoshi.5s@kyoto-u.jp}
}
\date{\today}

\begin{document}
\maketitle
\hrule
\begin{abstract}
We study tridiagonal pairs of type II. These involve two linear transformations $A$ and $A^\star$. We define two bases. In the first one, $A$ acts as a diagonal matrix while $A^\star$ acts as a block tridiagonal matrix, and in the second one, $A$ acts as a block tridiagonal matrix while $A^\star$ acts as a diagonal matrix. We obtain the change of basis coefficients between these two bases. The coefficients are special functions that are written as a nested product of polynomials that resemble Racah polynomials but involve shift operators in their expression.
\end{abstract}
\hrule

\medskip

\begin{center}
\begin{minipage}{16cm}

\textbf{MSC2020 database:} 47B36; 33C50; 33C80; 16G60
\end{minipage}
\end{center}

\vfill
\clearpage
\newpage



\section{Introduction}
\label{sec:intro}
\paragraph{Leonard pairs.}{}
Consider two finite-dimensional
matrices $A$ and $A^\star$.
An elementary result in linear algebra is:
\begin{center}
\textit{``Two diagonalizable matrices $A$ and $A^\star$ commute with each other if and only if they can be
diagonalized simultaneously.''}
\end{center}
The simplest departure from this trivial case is perhaps the situation where
upon diagonalizing $A$, the matrix $A^\star$ takes a tridiagonal form.
Now suppose that the opposite also holds---that is, $A$ takes a tridiagonal form
when diagonalizing $A^\star$---and in addition that both matrices are irreducible.
This is the setting of \textit{Leonard pairs}
\cite{Terwilliger1987, Terwilliger2001}.

In this looser but more general setting, the matrices $A$ and $A^\star$ do not
commute.
Instead, they obey well-studied algebraic relations known as the
\textit{Askey--Wilson relations} \cite{Zhedanov1991}:
\begin{subequations}\label{eq:AWrels}
\begin{align}
 & A^{2\phantom{\star}}A^{{\star}}
 -\beta A^{\phantom{\star}}A^{{\star}}A^{\phantom{\star}}
 +A^{{\star}}A^{2\phantom{\star}}
 -\gamma^{\phantom{\star}}\{A,A^\star\}
 -\rho^{\phantom{\star}} A^{{\star}}=
 \gamma^{\star}A^{2\phantom{\star}}
 +\omega A^{\phantom{\star}} + \eta^{\phantom{\star}}\,,\\
 & A^{{\star2}}A^{\phantom{\star}}
 -\beta A^{{\star}}A^{\phantom{\star}}A^{{\star}}
 +A^{\phantom{\star}}A^{{\star}2}
 -\gamma^\star\{A,A^\star\}
 -\rho^\star A^{\phantom{\star}}=
 \gamma^{\phantom{\star}}A^{{\star}2}
 +\omega A^{{\star}} + \eta^{{\star}}\,,
\end{align}
\end{subequations}
where $\beta$, $\gamma$, $\gamma^\star$, $\rho$, $\rho^\star$, $\omega$, $\eta$,
$\eta^\star$
are some parameters
(this has been shown by Terwilliger in \cite{Terwilliger1987}).
There are two natural bases in this situation:
the one where $A$ is diagonal, and the other one where $A^\star$ is diagonal.
These two natural bases are related by a change of basis.
The matrix that gives the change of basis is remarkable:
in the most general case \cite{Leonard1982},
its entries contain $q$-Racah polynomials, and various
limiting cases yield other families of the Askey scheme of hypergeometric
orthogonal polynomials \cite{KoekoekLeskyetal2010}.

This is a striking result.
From a simple linear algebra problem, one is led to the (finite)
orthogonal polynomials families of the Askey scheme.
And the matrices $A$ and $A^\star$
encode fundamental properties of these families:
their recurrence and difference equations.

\paragraph{Tridiagonal pairs.}
The generalization from commuting matrices to Leonard pairs was valuable.
It makes sense to look for further extensions.
One generalization of Leonard pairs is the concept of \textit{tridiagonal pairs}.
In this looser setting (see Definition \ref{def:TD}), one demands that
in the basis where $A$ is diagonal, $A^\star$ is block tridiagonal, and vice-versa.
The algebraic relations between $A$ and $A^\star$ have been obtained in the general
case \cite{Terwilliger1993a, Terwilliger2001a} and are called
\textit{tridiagonal relations}:
\begin{subequations}\label{eq:TDrels}
\begin{align}
 &\big[\,A^{\phantom{\star}},\
 A^{2\phantom{\star}}A^{{\star}}
 -\beta A^{\phantom{\star}}A^{{\star}}A^{\phantom{\star}}
 +A^{{\star}}A^{2\phantom{\star}}
 -\gamma^{\phantom{\star}}\{A,A^\star\}
 -\rho^{\phantom{\star}} A^{{\star}}\,
 \big]=0\,,\\
 &\big[\,A^\star,\
 A^{{\star2}}A^{\phantom{\star}}
 -\beta A^{{\star}}A^{\phantom{\star}}A^{{\star}}
 +A^{\phantom{\star}}A^{{\star}2}
 -\gamma^\star\{A,A^\star\}
 -\rho^\star A^{\phantom{\star}}\,
 \big]=0\,,
\end{align}
\end{subequations}
where $\beta$, $\gamma$, $\gamma^\star$, $\rho$, $\rho^\star$ are some parameters.
These relations, introduced for the study of the association scheme
\cite{Terwilliger1993a}, have also been studied in the context of integrable
systems \cite{BaseilhacKoizumi2005,Baseilhac2005a} and of the circular
Hessenberg pairs \cite{Lee22}. They also generalize the $q$-Dolan--Grady
relations \cite{Baseilhac2005} or the cubic $q$-Serre relations
\cite{ItoTerwilliger07}.
Inspired by the story of Leonard pairs, it would be natural to consider two bases
(one where $A$ is diagonal and another one where $A^\star$ is diagonal)
and look for the change of basis coefficients between the two bases.
In these should appear special functions which are generalizations of the
polynomials of the Askey scheme.

\paragraph{The goal.}
To the best of our knowledge, so far such results have not been obtained in general.
The question of what special functions arise in the study of tridiagonal pairs
still remains.
This is what we will clarify in this paper.
\begin{center}
\underline{
We will explicitly provide these special functions for tridiagonal
pairs of type II.
}
\end{center}
The precise definition of tridiagonal pairs of type II is given in
Definition \ref{def:typeII}.
Roughly speaking they correspond to a $q\to1$ limit of the most general case
of tridiagonal pairs.

\paragraph{The bases.}
There are a number of challenges to face in order to solve this change of basis
problem.
One is the fact that there are multiple choices of bases
where $A$ is diagonal and $A^\star$ is block tridiagonal (and vice-versa).
This is due to the degeneracy of the eigenvalues of $A$ and $A^\star$.
Thus, from the start, one cannot speak of ``two natural bases''.
Fortunately, there exists another basis called the \textit{split basis}.
The split basis, given in equation \eqref{eq:AAsVn}, corresponds to a basis
where $A$ takes an upper triangular, block bidiagonal form
while $A^\star$ takes a lower triangular, block bidiagonal form
(or upper bidiagonal and lower bidiagonal form in the case of Leonard pairs).
This basis was used to great effect by Terwilliger and proved to be a key
ingredient in his study of Leonard pairs.

The split basis provides a well-defined basis and we use it as our starting point.
From there, we obtain two natural bases.
One basis is such that $A$ is diagonal, $A^\star$ is block tridiagonal
and the change of basis from the split basis to this basis is given by a lower
triangular matrix with $1$s on the diagonal.
The other basis is obtained similarly by a change of basis involving an upper
triangular matrix with $1$s on the diagonal and leads to $A^\star$ being
diagonal while $A$ is block tridiagonal.
Equipped with two well-defined bases, the problem is now well posed and we solve it
in this paper.

\paragraph{Multivariate special functions from univariate ones.}
The special functions that are obtained have quite an intricate expression
(see equation \eqref{eq:cob_prod}) but upon introducing shift operators they can
be reworked into familiar functions, namely entangled products of polynomials
that closely resemble Racah polynomials (see equation \eqref{eq:cob_fact}).
They belong to the class of multivariate special functions, which have been
studied more and more in recent years
\cite{Griffiths1971, Tratnik1991, Rosengren1998, GasperRahman2005, Scarabotti2007, Grunbaum2007,
HoareRahman2008, GeronimoIliev2010, GeronimoIliev2011, Iliev2011, IlievTerwilliger2012, Iliev2012,
GenestVinetetal2013b, DiaconisGriffiths2014, Griffiths2016, GenestPostetal2017}.
These generalizations attracted even more attentions recently with
the introduction of multivariate $P$- and $Q$-polynomial association schemes
\cite{BernardCrampeetal2023, BannaiKuriharaetal2023}.

The problem of defining multivariate generalizations of the families of
the Askey scheme is a long standing one.
Various paths to define multivariate special functions have been explored.
The most popular approach is to attempt to write multivariate polynomials as sums
of entangled products of univariate polynomials.
Thus, expressions of the form \eqref{eq:cob_fact} look familiar, but the peculiarity
here is the presence of shift operators in the arguments of the generalized
hypergeometric functions.
We also highlight the work of Baseilhac, Vinet and Zhedanov
\cite{BaseilhacVinetetal2017} where entangled products of $q$-Racah polynomials
were identified as overlaps between some bases of the $q$-Onsager algebra
(note that no shift operators appear there).

Using univariate polynomials as Lego blocks with which to build multivariate
special functions through appropriate entangled products is a research direction
of increasing importance.
Recently, new methods to study and build these type of functions
have been designed \cite{CrampeFrappatetal2022,CrampeFrappatetal2023, CrampeZaimi2023,
CrampeFrappatetal2024, CrampeFrappatetal2024a}.
This is a topic that deserves more exploration and this paper adds important
insights.

\paragraph{Structure of the paper.}
In Section \ref{sec:tdpairs}, we define tridiagonal pairs of type II, which are
the mise-en-scène for this paper.
We introduce families of tridiagonal pairs of type II in the split basis
in Section \ref{sec:split} and highlight some of their properties
(involution, irreducibility, tridiagonal relations).
In Section \ref{sec:diag} we diagonalize $A$ and $A^\star$ by providing the
explicit change of basis from the split basis to the two bases.
The special functions arising as change of basis coefficients between the two
natural bases are obtained in Section \ref{sec:cob}.
Limiting cases are given in Section \ref{sec:limit}.
We conclude with some outlooks in Section \ref{sec:outlooks}.
The particular notations used throughout the paper are jotted down in Appendix
\ref{sec:notations}.

%
%


\section{Tridiagonal pairs of type II}
\label{sec:tdpairs}
We first recall the definition and main properties of tridiagonal pairs.
\begin{defi}[\cite{ItoTanabeetal2001}]
\label{def:TD}
Let $\mathbb{F}$ denote a field, and let $\cV$ denote a vector space
over $\mathbb{F}$ with finite positive dimension.
By a \emph{tridiagonal pair} (or TD-pair) on $\cV$, we mean
an ordered pair $(A, A^\star)$, where $A$ and $A^\star$ are linear
transformations on $\cV$ that satisfy the following four conditions:
\begin{itemize}
\item[(i)] $A$ and $A^\star$ are both diagonalizable on $\cV$.
\item[(ii)] There exists an ordering $\cV_0, \cV_1,\dots, \cV_d$ of the
eigenspaces of $A$ such that
\begin{align}
 A^\star\, \cV_i \subseteq \cV_{i-1}+\cV_i+\cV_{i+1}\,
 \qquad (0\leq i\leq d)
\end{align}
where $\cV_{-1}= 0,\ \cV_{d+1}= 0$.
\item[(iii)] There exists an ordering
$\cV^\star_0, \cV^\star_1,\dots, \cV^\star_\delta$ of the eigenspaces of
$A^\star$ such that
\begin{align}
 A\, \cV^\star_i \subseteq \cV^\star_{i-1}
 +\cV^\star_i+\cV^\star_{i+1}\, \qquad (0\leq i\leq \delta)
\end{align}
where $\cV^\star_{-1}= 0,\ \cV^\star_{\delta+1}= 0$.
\item[(iv)] There is no subspace $\cW$ of $\cV$ such that both
$A\cW\subseteq \cW$, $A^\star \cW\subseteq\cW$,
other than $\cW=0$ and $\cW=\cV$ (irreducibility).
\end{itemize}
\end{defi}
In this paper, we restrict ourselves to the case when the field $\mathbb{F}$ is
the field of complex numbers $\mathbb{C}$.
For a TD-pair, it is well-known that $A$ and $A^\star$ have the same number of
eigenspaces \textit{i.e.}~$d+1=d^\star+1$ ($d$ is called the \emph{diameter} of
the TD-pair) and that the following conditions hold
\cite{ItoTanabeetal2001}:
\begin{align}
 \rho_i=\text{dim}\cV_i = \text{dim}\cV^\star_i=\text{dim}\cV_{d-i}
       =\text{dim}\cV^\star_{d-i}\qquad (0\leq i \leq d)\,.
\end{align}
The tuple $(\rho_0,\rho_1,\dots,\rho_d)$ is called the \emph{shape} of a
TD-pair.
Because we restrict ourselves to a closed field (the complex numbers), the
TD-pair is sharp \textit{i.e.}~$\rho_0=1$.

It is conjectured that, for any TD-pair, there exist $N$ integers
$\ell_1,\ell_2,\dots,\ell_N\geq 0$, $N\geq 1$ such that the shape satisfies
the \emph{character formula}
\begin{align}\label{eq:rhol}
 \sum_{i=0}^{d} \rho_i \lambda^i
 = \prod_{p=1}^N \frac{1-\lambda^{\ell_p+1}}{1-\lambda}\,.
\end{align}
In particular, the diameter $d$ satisfies $d=\ell_1+\dots+ \ell_N$.
The case $N=1$ and $\Bell=(\ell)$ corresponds to
$\rho_0=\rho_1=\dots=\rho_{\ell}=1$
\textit{i.e.}~the spectra of $A$ and $A^\star$ are simple.
Then, a TD-pair becomes a Leonard pair.

\textbf{In this paper, we provide and study families of TD-pairs
for any choice of the $N$-tuple}
\begin{align}
 \Bell=(\ell_1,\dots,\ell_N)\,.
\end{align}
It is well-established that the eigenvalues $\theta_i$ and $\theta_i^\star$
of $A$ and $A^\star$ defining a TD-pair can take only few particular forms.
In this paper, we focus on the TD-pairs of type II.
\begin{defi}\label{def:typeII}
A tridiagonal pair (see Definition \ref{def:TD}) such that
the eigenvalues $\theta_i$, $\theta^\star_i$ of $A$, $A^\star$
can be expressed as
\begin{subequations}
\begin{align}
 \theta_i&=c_0 + c_1 i + c_2 i^2\,,\\
 \theta^\star_i&=c^\star_0 + c^\star_1 i + c^\star_2 i^2\,,
\end{align}
\end{subequations}
for some scalars $c_0,c_1,c_2,c^\star_0,c^\star_1,c^\star_2\in\mathbb{C}$
is called a \emph{tridiagonal pair of type II}.
\end{defi}

\section{Split bases for tridiagonal pairs of type II}
\label{sec:split}
In this section, we introduce a family of tridiagonal pairs of type II in
the so-called split basis.
The terminology \textit{split basis} refers to a basis in which
the two linear transformations $A$ and $A^\star$ act
as lower (resp. upper) block tridiagonal matrices.
We note that these families of TD-pairs also appeared in \cite{Ito2014}.
For $N=2$ and $\Bell=(\ell_1,1)$, we also recover the TD-pairs obtained in
\cite{ItoSato2014} (see Remark \ref{rem:IS}).

Below we provide the explicit actions of $A$ and $A^\star$ in the split basis.
We then prove that $A$ and $A^\star$ obey algebraic relations, the
tridiagonal relations.

\subsection{Split basis}
Now and for the rest of the paper, we fix a $N$-tuple
$\Bell=(\ell_1,\dots,\ell_N)$.
Let us denote by $V^{\bn}=V^{(n_1,\dots,n_N)}$, $0\leq n_p \leq \ell_p$,
the canonical basis vectors of
$\cV=\mathbb{C}^{\ell_1+1}\otimes\dots \otimes \mathbb{C}^{\ell_N+1}$
and take these vectors to be orthonormal
\begin{align}
 \langle V^{\mathbf{m}},V^{\bn} \rangle
 =\delta_{m_1,n_1}\dots \delta_{m_N,n_N}\,,
\end{align}
where $\langle V^{\mathbf{m}},V^{\bn} \rangle$ is the standard dot product.
By convention, $V^{\bn}=0$ as soon as there exists $1\leq i\leq N$
such that $n_i>\ell_i$ or $n_i<0$. We give in Appendix \ref{sec:notations}
the conventions used in this paper for the $N$-tuples.

For $i=0,1,\dots,|\Bell|$, write the \emph{eigenvalues} of type II
as follows:
\begin{align}\label{eq:thetadef}
 \theta_i=\theta_0 + hi(i+\omega)\,,\qqquad
 \theta^\star_i=\theta^\star_0 + h^\star i(i+\omega^\star)\,,
\end{align}
for some scalars
$\theta_0,\theta_0^\star,h,h^\star,\omega,\omega^\star\in \mathbb{C}$.
For $p=1,\dots,N$, we also define
\begin{align}
\begin{aligned}
 \xi_{\bn,p}^{\phantom{\star}\sss(N)} &= h^{\phantom{\star}}\,
 (|\bn|_1^{p-1}+|\bn|_1^p+|\Bell|_{p\phantom{+1}}^N+a_p
  +\omega^{\phantom{\star}})\, n_p\,, \\
 \xi^{\star \sss(N)}_{\bn,p} &= h^\star\,
 (|\bn|_1^{p-1}+|\bn|_1^p+|\Bell|_{p+1}^N-a_p
  +\omega^\star)\, (n_p-\ell_p)\,,
\end{aligned}
\end{align}
where $a_p\in \mathbb{C}$ are free parameters.

Let $A$, $A^\star$ be linear operators on $\cV$
whose action on the basis vectors is given by
\begin{subequations}\label{eq:AAsVn}
\begin{align}\label{eq:AAsVn1}
 A\, V^{\bn}
 &=\theta_{|\bn|} V^{\bn}
 +\sum_{p=1}^N \xi_{\bn+\be_p,p}^{\sss(N)}\
 V^{\bn+\be_p}\,,\\
 \label{eq:AAsVn2}
 A^\star\, V^{\bn}
  &=\theta^{\star}_{|\bn|}V^{\bn}
 +\sum_{p=1}^N \xi^{\star (N)}_{\bn-\be_p,p}\
 V^{\bn-\be_p}\,.
\end{align}
\end{subequations}
Due to the triangular form of the matrices $A$ and $A^\star$ in
the basis $V^{\bn}$, it is immediate that their eigenvalues are
respectively $\theta_{|\bn|}$ and $\theta^\star_{|\bn|}$.
The dimensions of the corresponding eigenspaces
$\rho_{|\bn|}$ satisfy relation \eqref{eq:rhol}.
\begin{exam}\label{rem:IS}
Consider the case $N=2$ and $\Bell=(\ell,1)$ which provides a tridiagonal pair
with shape $(1,\underbrace{2,2,\dots,2}_{\ell},1)$.
The actions of $A$ and $A^\star$ become
\vspace{-1em}
\begin{subequations}
\begin{align}
AV^{(n,0)}&=\theta_n V^{(n,0)} +\xi^{(2)}_{(n+1,0),1}V^{(n+1,0)}
+\xi^{(2)}_{(n,1),2}V^{(n,1)}\,,\\
AV^{(n-1,1)}&=\theta_n V^{(n-1,1)} +\xi^{(2)}_{(n,1),1}V^{(n,1)}\,,\\
A^\star V^{(n,1)}&=\theta^\star_{n+1} V^{(n,1)}
+\xi^{\star(2)}_{(n-1,1),1}V^{(n-1,1)}
+\xi^{\star(2)}_{(n,0),2}V^{(n,0)}\,,\\
A^\star V^{(n+1,0)}&=\theta^\star_{n+1} V^{(n+1,0)}
+\xi^{\star(2)}_{(n,0),1}V^{(n,0)}\,.
\end{align}
\end{subequations}
In order to recover relations (56)--(64) in \cite{ItoSato2014},
our parameters must be transformed as follows:
\begin{align}
&\omega\mapsto \omega+1\,,&&\qquad\omega^\star\mapsto \omega^\star+1\,,
&&\quad a_1\mapsto a+\frac{\omega^\star-\omega-\ell-1}{2}\,,
&&\quad a_2\mapsto b+\frac{\omega^\star-\omega}{2}-1\,.
\end{align}
In \cite{ItoSato2014}, it was also proven that there always exists such a split basis for all tridiagonal pairs of type II with shape $(1,\underbrace{2,2,\dots,2}_{\ell},1)$. For other shapes, similar results remain an open problem.
\end{exam}

\subsection{The involution \texorpdfstring{$S$}{S}}
The linear transformations $A$ and $A^\star$ are related to one another through
a simple transformation.
Introduce the operator $S$ which acts as follows:
\begin{align} \label{eq:defS}
 S\, V^{\bn}=V^{\Bell-\bn}\,.
\end{align}
The operator $S$ is an involution ($S^{-1}=S$).
It relates the linear transformations $A$ and $A^\star$, up to an appropriate
change of parameters:
\begin{subequations}\label{eq:SAS}
\begin{align}
&SAS = A^\star\Big|_{
\theta_0^\star \to \theta_0+h|\Bell|(|\Bell|+\omega),~~
h^\star \to h,~~
\omega^\star \to -\omega-2|\Bell|}\,,\\
&SA^\star S = A\Big|_{
\theta_0 \to \theta_0^\star+h^\star|\Bell|(|\Bell|+\omega^\star),~~
h \to h^\star,~~
\omega \to -\omega^\star-2|\Bell|}\,.
\end{align}
\end{subequations}
This simple transformation will be used to simplify some proofs.

\subsection{Parameter constraints and irreducibility}
In \cite{Ito2014}, it has been proven that \eqref{eq:AAsVn} is an instance of a
tridiagonal pair\footnote{
  To recover the notation used in \cite{Ito2014}
  from our notation, perform the following modifications on the parameters:
  \begin{gather*}
  h\mapsto 1,\qquad h^\star\mapsto 1,\qquad
  \theta_0\mapsto \tfrac{1}{4}(s+t-1-|\Bell|)(s+t+1-|\Bell|),\qquad
  \theta_0^\star\mapsto \tfrac{1}{4}(s-t-1-|\Bell|)(s-t+1-|\Bell|)\,,\\
  N \mapsto n,\qquad \omega \mapsto s+t-|\Bell|,\qquad
  \omega^\star \mapsto s-t-|\Bell|,\qquad
  a_i\mapsto a_i-\tfrac{1}{2}(1+\ell_i+\omega-\omega^\star)\,.
  \end{gather*}
  Our parameters $h,h^\star,\theta_0,\theta_0^\star$ are mapped to
  particular values in \cite{Ito2014} since in that paper
  they use a standardization
  (affine transformation of $A$ and $A^\star$).
}
if the parameters satisfy various constraints.

Firstly, the eigenvalues must satisfy
$\theta_i\neq\theta_j$ and $\theta^\star_i\neq\theta^\star_j$
if $i\neq j$ which leads to
\begin{align}\label{eq:cond1}
 h,h^\star\neq 0\,,\qquad
 \omega,\,\omega^\star \notin \{ -2|\Bell|+1, \dots, -2, -1 \}\,.
\end{align}
Secondly, the parameters $a_i$ satisfy the following set of constraints
\cite{Ito2014}, for $i=1,2,\dots,N$:
\begin{align}\label{eq:cond2}
a_i,\,a_i+\omega-\omega^\star,a_i-|\Bell|-\omega^\star,\,a_i+|\Bell|+\omega
\notin\{ -\ell_i,-\ell_i+1,\dots,-1\}\,.
\end{align}
Let us briefly recall the origin of those constraints.
These come from assuming that the off-diagonal entries of the matrices $A$,
$A^\star$ are never vanishing
(otherwise, invariant subspaces of $\cV$ could be constructed).
For Leonard pairs~\cite{Terwilliger2003} those constraints correspond to
$\varphi_i \neq 0$ and $\phi_i\neq 0$.

Finally, there is an additional set of constraints between the parameters $a_i$.
To describe them let us define the string $S^\pm(\ell,a)$ as
\begin{align}
S^\pm(\ell,a)=\left\{\pm\left(a+k+\tfrac{1}{2}(\omega-\omega^\star)\right)\
\Big|\ k=1,2,\dots,\ell\right\}\,.
\end{align}
Two strings $S^\epsilon(\ell,a)$ and $S^{\epsilon'}(\ell',a')$
(with $\epsilon,\epsilon'\in\{+,-\}$) are said to be
\textit{in general position} if either one is contained in the other or their
union is not a string.
The third constraint are then expressed as follows,
for any $i,j=1,2\dots,N$ and $\epsilon_i,\epsilon_j\in\{+,-\}$:
\begin{align}\label{eq:cond3}
S^{\epsilon_i}(\ell_i,a_i)\quad\text{and}\quad
S^{\epsilon_j}(\ell_j,a_j)\quad\text{are in general position}\,.
\end{align}
If the previous conditions \eqref{eq:cond1}, \eqref{eq:cond2} and
\eqref{eq:cond3} are satisfied, the matrices $A$ and $A^\star$ given by
\eqref{eq:AAsVn} are a tridiagonal pair.


We shall assume from now on that all the parameter constraints are obeyed
and therefore that the linear maps $A$, $A^\star$ provided in \eqref{eq:AAsVn}
are tridiagonal pairs.

\subsection{Tridiagonal relations}
Tridiagonal pairs obey algebraic relations known as the
\emph{tridiagonal relations}.
Below, we state these relations in our case \eqref{eq:AAsVn}
and provide a proof.
\begin{theo}
The operators $A$ and $A^\star$ given by \eqref{eq:AAsVn} satisfy the
tridiagonal relations:
\begin{subequations}\label{eq:TDrelations}
\begin{align}
 &\big[\,A^{\phantom{\star}},\
 A^{2\phantom{\star}}A^{{\star}}
 -\beta A^{\phantom{\star}}A^{{\star}}A^{\phantom{\star}}
 +A^{{\star}}A^{2\phantom{\star}}
 -\gamma^{\phantom{\star}}\{A,A^\star\}
 -\rho^{\phantom{\star}} A^{{\star}}\,
 \big]=0\,,
 \label{eq:TDrelation1}\\
 &\big[\,A^\star,\
 A^{{\star2}}A^{\phantom{\star}}
 -\beta A^{{\star}}A^{\phantom{\star}}A^{{\star}}
 +A^{\phantom{\star}}A^{{\star}2}
 -\gamma^\star\{A,A^\star\}
 -\rho^\star A^{\phantom{\star}}\,
 \big]=0\,,
 \label{eq:TDrelation2}
\end{align}
with
\begin{align}
\beta=2\,,\quad \gamma=2h\,,\quad \gamma^\star=2h^\star\,,\quad
\rho=h(h(\omega^{2}-1)-4\theta_0)\,,\quad
\rho^\star=h^\star(h^\star({\omega^\star}^2-1)-4\theta^\star_0)\,.
 \label{eq:TDrelation3}
\end{align}
\end{subequations}
Here $[X,Y]=XY-YX$ and $\{X,Y\}=XY+YX$.
\end{theo}
\proof
Introduce the raising and lowering operators, $R$ and $L$, defined as follows
\begin{align}\label{def:RL}
 R\, V^{\bn} =\sum_{p=1}^N  \xi^{(N)}_{\bn+\be_p,p}\
 V^{\bn+\be_p}\,,\qquad
 L\, V^{\bn} =\sum_{p=1}^N  \xi^{\star(N)}_{\bn-\be_p,p}\
 V^{\bn-\be_p}\,.
\end{align}
The operators $R$ and $L$ correspond to the off-diagonal parts of the operators
$A$ and $A^\star$ respectively.

To prove \eqref{eq:TDrelation1} with parameters \eqref{eq:TDrelation3}, it is
equivalent\footnote{
This equivalence is proven by direct computation, see also \cite{ItoSato2014}.
}
to prove that these operators satisfy
\begin{align}\label{eq:R3L}
 [R,[R,[R,L]]]V^{\bn}&=
 -6hh^\star(4|\bn|+\omega+\omega^\star+4)R^2V^{\bn}\,,
\end{align}
where we have rewritten the r.h.s~using the relation
\begin{align}
3(\theta_i\theta^\star_i -\theta_{i}\theta^\star_{i+1}
-\theta_{i+1}\theta^\star_{i} +\theta_{i+1}\theta^\star_{i+2}
+\theta_{i+2}\theta^\star_{i+1} -\theta_{i+2}\theta^\star_{i+2})
=-6hh^\star(4i+\omega+\omega^\star+4)\,.
\end{align}
We prove \eqref{eq:R3L} by induction on $N$.
For $N=1$, the tuple $\bn$ reduces to $\bn=(n)$ and one gets
by direct computation
\begin{align}
 [R,[R,[R,L]]]V^{(n)}&=\Big(
 \xi^{\star(1)}_{(n-1),1}\xi^{(1)}_{(n),1}\xi^{(1)}_{(n+1),1}\xi^{(1)}_{(n+2),1}
 -3\xi^{(1)}_{(n+1),1}\xi^{\star(1)}_{(n),1}\xi^{(1)}_{(n+1),1}
 \xi^{(1)}_{(n+2),1}\nonumber\\
 &\qquad
 +3\xi^{(1)}_{(n+1),1}\xi^{(1)}_{(n+2),1}\xi^{\star(1)}_{(n+1),1}
 \xi^{(1)}_{(n+2),1}
 -\xi^{(1)}_{(n+1),1}\xi^{(1)}_{(n+2),1}\xi^{(1)}_{(n+3),1}
 \xi^{\star(1)}_{(n+2),1}
 \Big)V^{(n+2)}\nonumber\\
 &=-6hh^\star(4n+\omega+\omega^\star+4)
 \xi^{(1)}_{(n+1),1}\xi^{(1)}_{(n+2),1}V^{(n+2)}\,.
\end{align}
We recognize the r.h.s.~of \eqref{eq:R3L} in the last line
which proves this relation for $N=1$.
Now suppose that \eqref{eq:R3L} holds for $N-1$.
Decompose the action of $R$ and $L$ as follows
\begin{subequations}
\begin{align}
 &R\,V^{\bn}=\underbrace{
 \sum_{p=1}^{N-1}
 \xi^{(N)}_{\bn+\be_p,p}\
 V^{\bn+\be_p}
 }_{:=R_{[N-1]}V^{\bn}}\
 +\ \underbrace{
 \xi^{(N)}_{\bn+\be_N,N}\  V^{\bn+\be_N}
 }_{:=R_{(N)} V^{\bn}}\,,\\
 &L\,V^{\bn}=\underbrace{
 \sum_{p=1}^{N-1}
 \xi^{\star(N)}_{\bn-\be_p,p}\
 V^{\bn-\be_p}
 }_{:=L_{[N-1]}V^{\bn}}\
 +\ \underbrace{
 \xi^{\star(N)}_{\bn-\be_N,N}\  V^{\bn-\be_N}
 }_{:=L_{(N)} V^{\bn}}\,.
\end{align}
\end{subequations}
Then, expand the l.h.s.~of \eqref{eq:R3L} and use the following properties:
\begin{subequations}
\begin{align}
 R_{(N)}R_{[N-1]} V^{\bn} &=
 \xi_{\bn+\be_N,N}^{(N)}\Big|_{\omega\to\omega+2}
 R_{[N-1]} V^{\bn+\be_N}\,,\\
 L_{(N)} R_{[N-1]} V^{\bn}&=
 \xi_{\bn-\be_N,N}^{\star(N)}\Big|_{\omega^\star\to\omega^\star+2}
 R_{[N-1]} V^{\bn-\be_N}\,,\\
 R_{(N)} L_{[N-1]} V^{\bn}&=
 \xi_{\bn+\be_N,N}^{(N)}\Big|_{\omega\to\omega-2}
 L_{[N-1]} V^{\bn+\be_N}\,.
\end{align}
\end{subequations}
Afterwards, $11$ types of terms remain.
We analyze them one by one.
This is a straightforward calculation using the expressions of
$\xi^{(N)}_{\bn,p}$ and $\xi^{\star(N)}_{\bn,p}$.
The first $10$ terms and their coefficients are presented in Table
\ref{tab:11terms}.
\begin{table}
\centering
\begin{tabular}{r|l}
Term & Coefficient in front \\
\hline\hline\\[-.5em]
$(R_{[N-1]})^3V^{\bn-\be_N}$ & $0$ \\[.5em]
$(R_{[N-1]})^2L_{[N-1]}V^{\bn+\be_N}$ & $0$ \\[.5em]
$R_{[N-1]}L_{[N-1]}R_{[N-1]}V^{\bn+\be_N}$ & $0$ \\[.5em]
$L_{[N-1]}(R_{[N-1]})^2V^{\bn+\be_N}$ & $0$ \\[.5em]
$R_{[N-1]}L_{[N-1]}V^{\bn+2\be_N}$ & $0$ \\[.5em]
$L_{[N-1]}R_{[N-1]}V^{\bn+2\be_N}$ & $0$ \\[.5em]
$L_{[N-1]}V^{\bn+3\be_N}$ & $0$ \\[.5em]
$V^{\bn+2\be_N}$ &
$-6hh^\star(4|\bn|_1^N+\omega+\omega^\star+4)
\xi_{\bn+\be_N,N}^{(N)}\xi_{\bn+2\be_N,N}^{(N)}$ \\[.5em]
$R_{[N-1]} V^{\bn+\be_N}$ &
$-6hh^\star(4|\bn|_1^N+\omega+\omega^\star+4)
 \left(\xi_{\bn+\be_N,N}^{(N)}
       +\xi_{\bn+2\be_N,N}^{(N)}\Big|_{\omega\to\omega+2}\right)$
\\[.5em]
$(R_{[N-1]})^2 V^{\bn}$ & $-6hh^\star(4n_N-2\ell_N)$ \\[.5em]
\hline
\end{tabular}
\caption{Terms resulting from the expansion of the l.h.s.~of \eqref{eq:R3L}}
\label{tab:11terms}
\end{table}
\par\noindent
Lastly, remarking that for $p=1,\dots,N-1$,
\begin{align}
\xi^{(N)}_{\bn,p}=\xi^{(N-1)}_{\bn,p}\Big|_{\omega\to \omega+\ell_N}
\qquad\text{and}\qquad \xi^{\star (N)}_{\bn,p}
=\xi^{\star(N-1)}_{\bn,p}\Big|_{\omega^\star\to \omega^\star+\ell_N}\,,
\end{align}
and using the recurrence hypothesis, the eleventh and last term
in the expansion is computed:
\begin{align}
 [R_{[N-1]},[R_{[N-1]},[R_{[N-1]},L_{[N-1]}]]]V^{\bn}&
 =-6hh^\star(4|\bn|_{1}^{N-1}+\omega+\omega^\star+2\ell_N+4)
 (R_{[N-1]})^2V^{\bn}\,.
\end{align}
\par\noindent
Putting together all the different factors appearing in the l.h.s.~of
\eqref{eq:R3L}, we obtain its r.h.s., which proves relation \eqref{eq:R3L}.

Relation \eqref{eq:TDrelation2} is proved by conjugating \eqref{eq:TDrelation1}
by the operator $S$ defined in \eqref{eq:defS} and
making use of relations \eqref{eq:SAS}.
\endproof


\section{Explicit diagonalization}
\label{sec:diag}
In this section, we show that $A$ and $A^\star$ are diagonalizable
(\textit{i.e.}~they satisfy condition (i) of Definition \ref{def:TD}) by
providing explicit transition matrices.

\subsection{Diagonalization of \texorpdfstring{$A$}{A}}
Let us consider the new vectors $V(\bx)$, for
$\bx=(x_1,\dots,x_N)$ an $N$-tuple of integers such that $0\leq x_i\leq \ell_i$
($i=1,\dots,N$), defined by
\begin{align}\label{eq:Vxm}
 V(\bx)&=\sum_{\bn=\bx}^{\Bell}
 \underbrace{
 \frac{(-1)^{|\bn|-|\bx|}}
 {(2|\bx|+\omega+1)_{|\bn|-|\bx|}}
 \prod_{p=1}^N \binom{n_p}{x_p}\
 (|\bn|_{1}^{p-1}\!+|\bx|_{1}^{p}+|\Bell|_{p}^{N}
  +a_p+\omega+1)_{n_p-x_p}
 }_{=:\sC^{\sss(N)}_{\bn,\bx}}\
 V^{\bn}\,.
\end{align}
In the previous formula, the sum must be understood pointwise
\textit{i.e.}~it is a sum over all $n_1,\dots,n_N$ satisfying
$x_p\leq n_p \leq \ell_p$ simultaneously for $p=1,\dots, N$:
\begin{align}
 \sum_{\bn=\bx}^{\Bell}=\sum_{n_1=x_1}^{\ell_1}\;\sum_{n_2=x_2}^{\ell_2}
 \cdots \sum_{n_N=x_N}^{\ell_N}\,.
\end{align}
We shall use this definition in the rest of the paper.
The condition \eqref{eq:cond1} on $\omega$ ensures that
$\sC^{\sss(N)}_{\bn,\bx}$ is always well-defined
since the denominator $(2|\bx|+\omega+1)_{|\bn|-|\bx|}$
does not vanish.
Let us also remark that $\sC^{\sss(N)}_{\bn,\bn}=1$.
Then, $\sC^{\sss(N)}$, the matrix with entries
$\sC^{\sss(N)}_{\bn,\bx}$,
can be seen as a lower triangular matrix of determinant $1$
(all $1$s on the diagonal)
that connects the bases $V(\bn)$ and $V^{\bn}$ as follows:
$V(\bn)=\sC^{\sss(N)}V^\bn$.
\begin{prop}\label{pro:AV}
The vectors $V(\bx)$ diagonalize $A$:
\begin{align}\label{eq:AVdiag}
 &A\, V(\bx) =\theta_{|\bx|}V(\bx) \,.
\end{align}
\end{prop}
\proof
Use the change of basis \eqref{eq:Vxm} to rewrite
the l.h.s.~of the above equation
and then use the action of $A$ on the basis vectors $V^{\bn}$
in \eqref{eq:AAsVn} to obtain:
\begin{align}
 A\,V(\bx)
 &=\sum_{\bn=\bx}^{\Bell}
 \sC^{(N)}_{\bn,\bx}
 \left[\theta_{|\bn|} V^{\bn}
 +\sum_{p=1}^N \xi_{\bn+\be_p,p}^{\scriptscriptstyle (N)}
 V^{\bn+\be_p}\right]\,.
\end{align}
For any $p=1,\dots,N$,
observe that $\sC^{\sss(N)}_{\bn,\bx}=0$ for any
$\bn=\bx-\be_p$
(this is due to the presence of the binomial coefficients
in the expression of $\sC^{\sss(N)}_{\bn,\bx}$)
and recall that $V^{\Bell+\be_p}=0$.
Thus, one can rewrite the above as
\begin{align}
 A\,V(\bx)
 &=\sum_{\bn=\bx}^{\Bell}
 \left[\theta_{|\bn|}\, \sC^{\sss(N)}_{\bn,\bx}
 +\sum_{p=1}^N \xi_{\bn,p}^{\sss(N)}\,
 \sC^{\sss(N)}_{\bn-\be_p,\bx}
 \right]V^{\bn}\,.
\end{align}
Noticing that $\theta_{|\bx|}-\theta_{|\bn|}
=h(|\bn|+|\bx|+\omega)(|\bx|-|\bn|)$,
it follows that equation \eqref{eq:AVdiag} holds if
\begin{align}\label{eq:conditionAV}
 \sum_{p=1}^N \xi_{\bn,p}^{\sss(N)}
 \frac{\sC^{\sss(N)}_{\bn-\be_p,\bx}}
      {\sC^{\sss(N)}_{\bn,\bx}}
 =\theta_{|\bx|}-\theta_{|\bn|}
 =h(|\bn|+|\bx|+\omega)(|\bx|-|\bn|)\,.
\end{align}
The latter equation is proven by induction on the variable $N$.
First, one checks by direct computation that \eqref{eq:conditionAV} holds
for $N=1$.
Then, suppose that \eqref{eq:conditionAV} holds for a given $N-1$
and study the l.h.s.~of \eqref{eq:conditionAV}:
\begin{align}
 \sum_{p=1}^{N} \xi_{\bn,p}^{\sss(N)}
 \frac{\sC^{\sss(N)}_{\bn-\be_p,\bx}}
      {\sC^{\sss(N)}_{\bn,\bx}}
 =\underbrace{
 \sum_{p=1}^{N-1} \xi_{\bn,p}^{\sss(N)}
 \frac{\sC^{\sss(N)}_{\bn-\be_p,\bx}}
      {\sC^{\sss(N)}_{\bn,\bx}}}_{=:\cA}
 +\underbrace{
 \xi_{\bn,p}^{\sss(N)}
 \frac{\sC^{\sss(N)}_{\bn-\be_{N},\bx}}
      {\sC^{\sss(N)}_{\bn,\bx}}}_{=:\cB}\,.
\end{align}
By using the expression for $\sC^{\sss(N-1)}_{\bn,\bx}$,
observe that
\begin{align}\label{eq:term_par}
&\begin{aligned}
 \cA&=\frac{|\bn|_1^{N}+|\bx|_1^{N}+\omega}
                  {|\bn|_1^{N-1}+|\bx|_1^{N-1}+\omega}\
 \frac{|\bn|_1^{N-1}+|\bx|_1^{N}+\ell_{N}+a_{N}+\omega}
      {|\bn|_1^{N}+|\bx|_1^{N-1}+\ell_{N}+a_{N}+\omega}
 \left.\left(
 \sum_{p=1}^{N-1} \xi_{\bn,p}^{\sss(N-1)}
 \frac{\sC^{\sss(N-1)}_{\bn-\be_p,\bx}}
      {\sC^{\sss(N-1)}_{\bn,\bx}}
 \right)\right|_{\ell_{N-1}\mapsto\ell_{N-1}+\ell_{N}}\,\\
 &=h\frac{|\bn|_1^{N}+|\bx|_1^{N}+\omega}
    {|\bn|_1^{N}+|\bx|_1^{N-1}+\ell_{N}+a_{N}+\omega}
    (|\bx|_1^{N-1}-|\bn|_1^{N-1})
    (|\bn|_1^{N-1}+|\bx|_1^{N}+\ell_{N}+a_{N}+\omega)\,,
\end{aligned}\\
&\begin{aligned}
\cB&=h\frac{|\bn|_1^{N}+|\bx|_1^{N}+\omega}
    {|\bn|_1^{N}+|\bx|_1^{N-1}+\ell_{N}+a_{N}+\omega}
    (x_{N}-n_{N})
    (|\bn|_1^{N-1}+|\bn|_1^{N}+\ell_{N}+a_{N}+\omega)\,.
\end{aligned}
\end{align}
where the term in parenthesis in \eqref{eq:term_par} is evaluated
using \eqref{eq:conditionAV}.
Adding up both $\cA$ and $\cB$ and using
\begin{align}
 \frac{(x-m)(x+y+m+c)+(x+y-m-n)(2m+n+c)}{x+m+n+c}=x+y-m-n\,,
\end{align}
one arrives at
\begin{align}
 \sum_{p=1}^{N} \xi^{\sss(N)}_{\bn,p}
 \frac{\sC^{\sss(N)}_{\bn-\be_p,\bx}}
      {\sC^{\sss(N)}_{\bn,\bx}}
 =h(|\bn|_1^{N}+|\bx|_1^{N}+\omega)
  (|\bx|_1^{N}-|\bn|_1^{N})\,
\end{align}
which completes the proof.
\endproof

\begin{prop}[]
Recall the change of basis presented in \eqref{eq:Vxm}.
The inverse change of basis is given by
\begin{align}\label{eq:VmxJ}
 V^{\bn}&=\sum_{\bx=\bn}^{\Bell}
\underbrace{
 \frac{1}{(|\bn|+|\bx|+\omega)_{|\bx|-|\bn|}}
 \prod_{p=1}^N \binom{x_p}{n_p}
 (|\bn|_{1}^{p}+|\bx|_{1}^{p-1}+|\Bell|_{p}^{N}
  +a_p+\omega+1)_{x_p-n_p}}_{=:\overline{\sC}^{\sss(N)}_{\bx,\bn}}
 V(\bx)\,.
\end{align}
\end{prop}
\proof
Recall ${\sC}^{\sss(N)}_{\bn,\bn}=1$, that is, the determinant of
${\sC}^{\sss(N)}$ is equal to $1$.
Since $\overline{\sC}^{\sss(N)}_{\bn,\bn}=1$
(the determinant of $\overline{\sC}^{\sss(N)}$ is also equal to $1$),
to prove the proposition it remains to
show that the action of $A$ on both sides of \eqref{eq:VmxJ} yields
the expected formula \eqref{eq:AAsVn} \textit{i.e.}~we have to prove that
\begin{align}
 \theta_{|\bn|} \sum_{\bx=\bn}^{\Bell}\overline{\sC}^{\sss(N)}_{\bx,\bn}V(\bx)
 +\sum_{p=1}^N \xi_{\bn+\be_p,p}^{\sss(N)}\sum_{\bx=\bn+\be_p}^{\Bell}
 \overline{\sC}^{\sss(N)}_{\bx,\bn+\be_p}V(\bx)
 =\sum_{\bx=\bn}^{\Bell}\overline{\sC}^{\sss(N)}_{\bx,\bn} \theta_{|\bx|} V(\bx)\,.
\end{align}
Owing to the fact that $\overline{\sC}^{\sss(N)}_{\bn,\bn+\be_p}V(\bx)=0$,
rewrite the sum $\sum_{\bx=\bn+\be_p}^{\Bell}$ as $\sum_{\bx=\bn}^{\Bell}$.
Comparing the coefficients in front of each $V(\bx)$ on both sides,
it follows that the proposition is proven if
\begin{align}\label{eq:conditionAVn}
 \sum_{p=1}^N \xi_{\bn+\be_p,p}^{\sss(N)}
 \frac{\overline{\sC}^{\sss(N)}_{\bx,\bn+\be_p}}
      {\overline{\sC}^{\sss(N)}_{\bx,\bn}}
 =\theta_{|\bx|}-\theta_{|\bn|}\,.
\end{align}
Fortunately, this follows from the previous proposition.
Indeed, each term in the sum on the l.h.s.~can be rewritten as
\begin{align}
\begin{aligned}
 \xi_{\bn+\be_p,p}^{\sss(N)}
 \frac{\overline{\sC}^{\sss(N)}_{\bx,\bn+\be_p}}
      {\overline{\sC}^{\sss(N)}_{\bx,\bn}}
 =h(|\bn|+|\bx|+\omega)(x_p-n_p)
 &\frac{|\bn|_1^{p-1}+|\bn|_1^{p}+|\Bell|_p^{N}+a_p+\omega+1}
      {|\bn|_1^{p}+|\bx|_1^{p-1}+|\Bell|_p^{N}+a_p+\omega+1}\\
\times \prod_{k=p+1}^{N}
 &\frac{|\bn|_1^{k-1}+|\bx|_1^{k}+|\Bell|_k^{N}+a_k+\omega+1}
      {|\bn|_1^{k}+|\bx|_1^{k-1}+|\Bell|_k^{N}+a_k+\omega+1}\,,
\end{aligned}
\end{align}
which, upon closer examination, is precisely equal to
\begin{align}
 \xi_{\bn+\be_p,p}^{\sss(N)}
 \frac{\overline{\sC}^{\sss(N)}_{\bx,\bn+\be_p}}
      {\overline{\sC}^{\sss(N)}_{\bx,\bn}}
 =\left.\left(
 \xi_{\bn,p}^{\sss(N)}
 \frac{{\sC}^{\sss(N)}_{\bn-\be_p,\bx}}
      {{\sC}^{\sss(N)}_{\bn,\bx}}
  \right)\right|_{\ell_N\mapsto\ell_N+1}\,.
\end{align}
Thus, using \eqref{eq:conditionAV} which was proved in the previous proposition,
the result follows.
\endproof

\begin{rema}
For $0\leq i \leq d=|\Bell|$, the eigenspace $\cV_i$ satisfying
the condition (ii) of Definition \ref{def:TD} is
\begin{align}
 \cV_i= \text{span}_{\mathbb{C}}
 \big\{ V(\bx) \ \big|\ |\bx|=i\big\}\,.
\end{align}
Therefore, the action of $A^\star$ on these eigenspaces is block tridiagonal.
Using relation \eqref{eq:Vxm},
then using the action \eqref{eq:AAsVn2} of $A^\star$ on $V^{\bn}$
and finally expressing the vectors $V^{\bn}$ in terms of $V(\bx)$
using \eqref{eq:VmxJ}, one obtains
\begin{align}
A^\star\,V(\bx)=
\sum_{\bn=\bx}^{\Bell}\sum_{\by=\bn}^{\Bell}
\theta^\star_{|\bn|}
{\sC}^{\sss(N)}_{\bn,\bx} \overline{\sC}^{\sss(N)}_{\by,\bn}
V(\by)
+
\sum_{\bn=\bx}^{\Bell}\sum_{p=1}^N\sum_{\by=\bn-\be_p}^{\Bell}
\xi^{\star \sss(N)}_{\bn-\be_p,p}
{\sC}^{\sss(N)}_{\bn,\bx} \overline{\sC}^{\sss(N)}_{\by,\bn-\be_p}
V(\by)\,.
\end{align}
The block tridiagonal action can then be read
off the above expression:
\begin{align}\label{eq:AsBT}
\begin{aligned}
A^\star\,V(\bx)&=
\sum_{p=1}^N \xi^{\star \sss(N)}_{\bx-\be_p,p}V(\bx-\be_p)
+\sum_{p=1}^{N}\left(
\theta^\star_{|\bx|}\overline{\sC}^{\sss(N)}_{\bx+\be_p,\bx}
+\theta^\star_{|\bx+\be_p|}{\sC}^{\sss(N)}_{\bx+\be_p,\bx}
\right)V(\bx+\be_p)\\
&+\theta^\star_{|\bx|}V(\bx)+\sum_{p,q=1}^{N}\left(
\xi^{\star \sss(N)}_{\bx-\be_q,q}
\overline{\sC}^{\sss(N)}_{\bx+\be_p-\be_q,\bx-\be_q}
+\xi^{\star \sss(N)}_{\bx+\be_p-\be_q,q}{\sC}^{\sss(N)}_{\bx+\be_p,\bx}
\right)V(\bx+\be_p-\be_q)\\
&+\sum_{p=1}^N \sum_{1\leq q\leq r\leq N}
\sum_{\bn=\bx}^{\bx+\be_q+\be_r} \xi^{\star \sss(N)}_{\bn-\be_p,p}
{\sC}^{\sss(N)}_{\bn,\bx}
\overline{\sC}^{\sss(N)}_{\bx+\be_q+\be_r-\be_p,\bn-\be_p}
V(\bx+\be_q+\be_r-\be_p)\,.
\end{aligned}
\end{align}
It is understood that the coefficients in front of other vectors outside the
block tridiagonal zone, such as $V(\bx+\be_p+\be_q)$, are zero by virtue of
condition (ii) of Definition \ref{def:TD} for tridiagonal pairs.
\end{rema}

\subsection{Diagonalization of \texorpdfstring{$A^\star$}{A*}}
Let us consider the new vectors $V_{\mathbf{i}}$,
for $\mathbf{i}=(i_1,\dots, i_N)$ an $N$-tuple of integers such that
$0\leq i_k \leq \ell_k$ ($k=1,\dots,N$),
defined by
\begin{align}\label{eq:VimN}
V_{\mathbf{i}}
&=\sum_{\bn=\mathbf{0}}^{\mathbf{i}}
\underbrace{\frac{(-1)^{|\bn|-|\mathbf{i}|}}
           {(|\mathbf{i}|+|\bn|+\omega^\star)_{|\mathbf{i}|-|\bn|}}
\prod_{p=1}^N \binom{\ell_p-n_p}{\ell_p-i_p}
(|\mathbf{i}|_{1}^{p-1}+|\bn|_{1}^{p}+|\Bell|_{p+1}^N-a_p+\omega^\star)
_{i_p-n_p} }_{=:\sD^{\sss(N)}_{\bn,\mathbf{i}}}  V^{\bn}\,.
\end{align}
The condition \eqref{eq:cond1} on $\omega^\star$ ensures that the previous
relation is always well-defined since the denominator
$(|\mathbf{i}|+|\bn|+\omega^\star)_{|\mathbf{i}|-|\bn|}$ does not vanish.
\begin{prop}
The vectors $V_{\mathbf{i}}$ diagonalize $A^\star$:
\begin{align}
 &A^\star V_{\mathbf{i}} =\theta^\star_{|\mathbf{i}|}V_{\mathbf{i}} \,.
\end{align}
\end{prop}
\proof
The proof can be done following the same lines as Proposition \ref{pro:AV}.
However, using the map $S$, we can show that $V_{\mathbf{i}}$ is an eigenvector
of $A^\star$ if the coefficients $\sD^{(N)}_{\bn,\mathbf{i}}$ satisfy
\begin{equation}
 \sD^{\sss(N)}_{\bn,\mathbf{i}}=
 \sC^{\sss(N)}_{\Bell-\bn,\Bell-\mathbf{i}}
 \Big|_{\omega\to -\omega^\star-2|\Bell|}\,.
\end{equation}
Using the explicit form of $\sC^{\sss(N)}_{\bn,\bx}$ given in \eqref{eq:Vxm}, we
prove the proposition.
\endproof
\begin{rema}
For $0\leq j \leq d=|\Bell|$,  the eigenspace $\cV^\star_j$
satisfying the condition (iii) of Definition \ref{def:TD} is
\begin{align}
 \cV^\star_j= \text{span}_{\mathbb{C}} \{ V_\mathbf{i} \ |\ |\mathbf{i}|=j\}\,.
\end{align}
The action of $A$ on these eigenspaces is block tridiagonal and is given below.
It is obtained in the same way as what was done in the above subsection
or alternatively from \eqref{eq:AsBT} by making use of the involution $S$.
\begin{align}\label{eq:ABT}
\begin{aligned}
A\,V_\bi&=
\sum_{p=1}^N \xi^{\sss(N)}_{\bi+\be_p,p}V_{\bi+\be_p}
+\sum_{p=1}^{N}\left(
\theta_{|\bi|}\overline{\sD}^{\sss(N)}_{\bi-\be_p,\bi}
+\theta_{|\bi-\be_p|}{\sD}^{\sss(N)}_{\bi-\be_p,\bi}
\right)V_{\bi-\be_p}\\
&+\theta_{|\bi|}V_{\bi}+\sum_{p,q=1}^{N}\left(
\xi^{\sss(N)}_{\bi+\be_q,q}
\overline{\sD}^{\sss(N)}_{\bi-\be_p+\be_q,\bi+\be_q}
+\xi^{\sss(N)}_{\bi-\be_p+\be_q,q}{\sD}^{\sss(N)}_{\bi-\be_p,\bi}
\right)V_{\bi-\be_p+\be_q}\\
&+\sum_{p=1}^N \sum_{1\leq q\leq r\leq N}
\sum_{\bn=\bi-\be_q-\be_r}^{\bi} \xi^{\sss(N)}_{\bn+\be_p,p}
{\sD}^{\sss(N)}_{\bn,\bi}
\overline{\sD}^{\sss(N)}_{\bi+\be_p-\be_q-\be_r,\bn+\be_p}
V_{\bi+\be_p-\be_q-\be_r}\,.
\end{aligned}
\end{align}
\end{rema}
\begin{prop}
The inverse change of basis is given as follows:
\begin{align}
V^{\bn}
&=\sum_{\bi=\mathbf{0}}^{\bn}
\underbrace{\frac{1}
           {(2|\bi|+\omega^\star+1)_{|\bn|-|\bi|}}
\prod_{p=1}^N \binom{\ell_p-i_p}{\ell_p-n_p}
(|\bn|_{1}^{p-1}+|\bi|_{1}^{p}+|\Bell|_{p+1}^N-a_p+\omega^\star)
_{n_p-i_p} }_{=:\overline{\sD}^{\sss(N)}_{\bi,\bn}}  V_{\bi}\,. \label{eq:Vni}
\end{align}
\end{prop}
\proof
This is computed straightforwardly using the same trick as previous proposition.
\endproof

\section{Change of basis and special functions}
\label{sec:cob}
In the theory of Leonard pairs, two natural bases are those where each linear
map acts diagonally.
Orthogonal polynomials of the ($q$-)Askey scheme appear as entries of the change
of basis matrices between these two natural bases.
Similarly, here we are interested in the special functions that appear as change
of basis coefficients between two natural bases for tridiagonal pairs.

Concretely, we want to compute the elements $T_{\mathbf{i}}(\bx)$ of the change
of basis between $V(\bx)$ and $V_{\mathbf{i}}$:
\begin{align}
 V_{\mathbf{i}}=\sum_{\bx=\mathbf{0}}^{\Bell}
 T_{\mathbf{i}}(\bx)\, V(\bx)\,.
\end{align}
Remarking that relations \eqref{eq:VmxJ} and \eqref{eq:VimN} imply
\begin{align}
 V_{\mathbf{i}}
&=\sum_{\bn=\mathbf{0}}^{\mathbf{i}}\sD^{\sss(N)}_{\bn,\mathbf{i}}
\sum_{\bx=\bn}^{\Bell}
\overline{\sC}^{\sss(N)}_{\bx,\bn}
V(\bx)=
\sum_{\bx=\mathbf{0}}^{\Bell}\;
\sum_{\bn=\mathbf{0}}^{\text{min}(\mathbf{i},\bx)}
\sD^{\sss(N)}_{\bn,\mathbf{i}} \overline{\sC}^{\sss(N)}_{\bx,\bn}\;V(\bx)\,,
\end{align}
where $\text{min}(\mathbf{i},\bx)\equiv(\text{min}(i_1,x_1),\dots,
\text{min}(i_N,x_N))$,
one gets
\begin{align}\label{eq:TDC}
 T_{\mathbf{i}}(\bx)=\sum_{\bn=\mathbf{0}}^{\text{min}(\mathbf{i},\bx)}
 \sD^{\sss(N)}_{\bn,\mathbf{i}} \;\overline{\sC}^{\sss(N)}_{\bx,\bn}\,.
\end{align}
The previous expression of $T_{\mathbf{i}}(\bx)$ can be simplified and expressed
in terms of generalized hypergeometric functions.
For $i,x,\ell$ nonnegative integers satisfying $i,x\leq \ell$, we introduce the
function $R$ as follows:
\begin{align}\label{eq:4F3}
\rac{i,x,a_1,a_2}{\ell,b_1,b_2}{Z}
&=\binom{\ell}{i}(b_1)_i(b_2)_x\;
\pFq{4}{3}{-i,\;-x,\; a_1,\; a_2}{b_1,\; b_2,\;-\ell}{Z}\nonumber\\
 &=\binom{\ell}{i}(b_1)_i(b_2)_x\sum_{k=0}^{\text{min}(i,x)}
  \frac{(-i,-x,a_1,a_2)_k}{(1,b_1,b_2,-\ell)_k}\, Z^k\,.
\end{align}
As the $Z$ might be an operator that does not commute with the coefficients in
the sum (see below), we need to fix a convention.
Here we use the convention that $Z^k$ is on the far right of each term in
the sum.
We also introduce the shift operator $e^{\partial_x}$ acting as follows:
\begin{equation}
 e^{\partial_x}f(x)=f(x+1)\,,
\end{equation}
and the function identity: $\mathbf{1}:x\mapsto 1$ such that $e^{\partial_x}\mathbf{1}(x)=1$.
\begin{theo}[\textbf{Change of basis between $V(\bx)$ and $V_{\mathbf{i}}$}]
The elements $T_{\mathbf{i}}(\bx)$ can be expressed as the following product:
\begin{align}\label{eq:cob_fact}
\begin{aligned}
 T_{\mathbf{i}}(\bx)&=
 \frac{(-1)^{|\mathbf{i}|}}{(|\mathbf{i}|+\omega^\star)_{|\mathbf{i}|}
 (|\bx|+\omega)_{|\bx|}} \\
 &\times{\overrightarrow{\prod}}_{p=1}^{N}
 \rac{i_p,\;x_p,\; |\mathbf{i}|+\omega^\star,\; |\bx|+\omega}
 {\ell_p,\; |\mathbf{i}|_{1}^{p-1}+|\Bell|_{p+1}^{N}+\omega^\star-a_p,\;
  |\bx|_{1}^{p-1}+|\Bell|_{p}^{N}+\omega+a_p+1}
 {e^{\partial_{i_p}+\partial_{x_p}}}
\cdot\mathbf{1}(i_N)\mathbf{1}(x_N)\,.
\end{aligned}
\end{align}
where the product is ordered,
\textit{i.e.}~${\overrightarrow{\prod}}_{p=1}^{N}X_p=X_1X_2\dots X_{N}$.
\end{theo}
\proof
Using properties of the Pochhammer symbol, the coefficients
$\overline{\sC}^{\sss(N)}_{\bx,\bn}$ and $\sD^{\sss(N)}_{\bn,\mathbf{i}}$
can be rewritten as follows:
\begin{align}\label{eq:Cb}
\overline{\sC}^{\sss(N)}_{\bx,\bn}=
\frac{(-1)^{|\bn|}}{(|\bx|+\omega)_{|\bx|}}
\prod_{p=1}^N
\frac{(-x_p,\ |\bx|+|\bn|_{1}^{p-1}+\omega)_{n_p}
      (|\bn|_{1}^{p-1}+|\bx|_{1}^{p-1}+|\Bell|_{p}^{N}+a_p+\omega+1)_{x_p}}
     {(1,\ |\bn|_{1}^{p-1}+|\bx|_{1}^{p-1}+|\Bell|_{p}^{N}+a_p+\omega+1)_{n_p}}
\end{align}
and
\begin{align}\label{eq:Db}
\sD^{\sss(N)}_{\bn,\mathbf{i}}
=\frac{(-1)^{|\bn|-|\mathbf{i}|}}{(|\mathbf{i}|+\omega^\star)_{|\mathbf{i}|}}
\prod_{p=1}^N \binom{\ell_p}{i_p}
\frac{(-i_p,\ |\mathbf{i}|+|\bn|_{1}^{p-1}+\omega^\star)_{n_p}
      (|\mathbf{i}|_{1}^{p-1}+|\bn|_{1}^{p-1}+|\Bell|_{p+1}^N-a_p
       +\omega^\star)_{i_p}}
     {(-\ell_p,\ |\mathbf{i}|_{1}^{p-1}+|\bn|_{1}^{p-1}+|\Bell|_{p+1}^N-a_p
      +\omega^\star)_{n_p}}\,.
\end{align}
Putting these together, we can write $T_{\bi}(\bx)$ as follows:
\begin{align}\label{eq:cob_prod}
\begin{aligned}
T_{\bi}(\bx)=\sum_{\bn=\mathbf{0}}^{\min(\bi,\bx)}
\prod_{p=1}^{N}
& \
\frac{(-x_p)_{n_p}(-i_p)_{n_p}(-\ell_p)_{i_p}}
     {(1)_{n_p}(-\ell_p)_{n_p}(1)_{i_p}}
\frac{(|\bx|_1^{p-1}+|\bn|_1^p+|\Bell|_{p}^N+a_p+\omega+1)_{x_p-n_p}}
     {(|\bx|+|\bn|+|\bx|_1^{p-1}-|\bn|_1^{p-1}+\omega)_{x_p-n_p}}\\
&\quad\times
\frac{(|\bi|_1^{p-1}+|\bn|_1^p+|\Bell|_{p+1}^N-a_p+\omega^\star)_{i_p-n_p}}
     {(|\bi|+|\bn|+|\bi|_1^{p-1}-|\bn|_1^{p-1}+\omega^\star)_{i_p-n_p}}\,.
\end{aligned}
\end{align}
This expression can be reworked into a more familiar form using the following trick.
Using the shift operators $e^{\partial_{x_p}}$, all the
dependence on $n_p$ of the expression \eqref{eq:Cb} can be gathered in
$e^{n_p\partial_{x_p}}$ as shown in the following expression:
\begin{align}
 \overline{\sC}^{\sss(N)}_{\bx,\bn}=
 \frac{(-1)^{|\bn|}}{(|\bx|+\omega)_{|\bx|}}
 {\overrightarrow{\prod}}_{p=1}^{N}
 & \
 \frac{(-x_p,|\bx|+\omega)_{n_p}(|\bx|_{1}^{p-1}+|\Bell|_{p}^{N}
  +a_p+\omega+1)_{x_p}}{(1,|\bx|_{1}^{p-1}+|\Bell|_{p}^{N}
  +a_p+\omega+1)_{n_p}}e^{n_p\partial_{x_p}}\cdot\mathbf{1}(x_N)\,.
\end{align}
We recall that the product is ordered, as mentioned at the end of the theorem.
Similar expressions can be obtained for $\sD^{\sss(N)}_{\bn,\mathbf{i}}$ using
the shift operator $e^{\partial_{i_p}}$. The result of the theorem is then proven
inserting these expressions in \eqref{eq:TDC}.

\endproof
It is immediate to see that the case $N=1$ yields univariate Racah polynomials.

\begin{theo}[\textbf{Inverse change of basis}]
The change of basis between $V_{\mathbf{i}}$ and $V(\bx)$ is written as
\begin{align}
 V(\bx)=\sum_{\bi=\mathbf{0}}^{\Bell}U_{\bi}(\bx)\, V_{\bi}\,,
\end{align}
with the coefficients given by
\begin{align}\label{eq:cobi_prod}
\begin{aligned}
U_{\bi}(\bx)&=
{\overrightarrow{\prod}}_{p=1}^{N}
\frac{(-1)^{i_p+\ell_p}}{(2|\bx|+|\Bell|_1^{p-1}-|\bx|_1^{p-1}+\omega+1)_{\ell_p-x_p}(2|\bi|+|\Bell|_1^{p-1}-|\bi|_1^{p-1}+\omega^\star+1)_{\ell_p-i_p}}\\
&\times \rac{\ell_p-x_p,\;\ell_p-i_p,\; -2|\bx|-|\Bell|_1^p+|\bx|_1^p-\omega,\; -2|\bi|-|\Bell|_1^p+|\bi|_1^p-\omega^\star}
 {\ell_p,\;  -|\Bell|-|\bx|_1^{p-1}-\ell_p-a_p-\omega,\;
  -|\Bell|-|\bi|_1^{p-1}+a_p+1-\omega^\star}
 {e^{\partial_{i_p}+\partial_{x_p}}}\cdot\mathbf{1}(x_N)\mathbf{1}(i_N)\,.
\end{aligned}
\end{align}
\end{theo}
\proof Using \eqref{eq:Vxm} and \eqref{eq:Vni}, one gets the following expression
\begin{align}\label{eq:cobi_prodbis}
\begin{aligned}
U_{\bi}(\bx)=\sum_{\bn=\max(\bx,\bi)}^{\Bell}
\prod_{p=1}^{N}& \
\frac{(-n_p)_{x_p}(-n_p)_{i_p}(-\ell_p)_{n_p}}{(1)_{x_p}(-\ell_p)_{i_p}(1)_{n_p}}
\frac{(|\bx|_1^p+|\bn|_1^{p-1}+|\Bell|_p^N+a_p+\omega+1)_{n_p-x_p}}
     {(2|\bx|+|\bn|_1^{p-1}-|\bx|_1^{p-1}+\omega+1)_{n_p-x_p}}\\
&\quad\times
\frac{(|\bi|_1^p+|\bn|_1^{p-1}+|\Bell|_{p+1}^N-a_p+\omega^\star)_{n_p-i_p}}
     {(2|\bi|+|\bn|_1^{p-1}-|\bi|_1^{p-1}+\omega^\star+1)_{n_p-i_p}}\,.
\end{aligned}
\end{align}
Performing the following change $\bn\to \Bell-\bn$ in the previous sum, one gets
\begin{align}\label{eq:cobi_prod2}
\begin{aligned}
U_{\bi}(\bx)=\sum_{\bn=\mathbf{0}}^{\min(\Bell-\bx,\Bell-\bi)}
\prod_{p=1}^{N}\ &
\frac{(n_p-\ell_p)_{x_p}(|\bx|_1^p-|\bn|_1^{p-1}+|\Bell|+a_p+\omega+1)_{\ell_p-n_p-x_p}}
{(1)_{x_p}(2|\bx|+|\Bell|_1^{p-1}-|\bn|_1^{p-1}-|\bx|_1^{p-1}+\omega+1)_{\ell_p-n_p-x_p}}
\\
&\times
\frac{(n_p-\ell_p)_{i_p}(-\ell_p)_{\ell_p-n_p}(|\bi|_1^p-|\bn|_1^{p-1}-\ell_p+|\Bell|-a_p+\omega^\star)_{\ell_p-n_p-i_p}}
{(-\ell_p)_{i_p}(1)_{\ell_p-n_p}(2|\bi|-|\bn|_1^{p-1}+|\Bell|_1^{p-1}-|\bi|_1^{p-1}+\omega^\star+1)_{\ell_p-n_p-i_p}}\,.
\end{aligned}
\end{align}
Then using identities for Pochhammer symbols one obtains
\begin{align}\label{eq:cobi_prod3}
\begin{aligned}
U_{\bi}(\bx)&=\sum_{\bn=\mathbf{0}}^{\min(\Bell-\bx,\Bell-\bi)}
\prod_{p=1}^{N}
\frac{(-\ell_p)_{x_p}(x_p-\ell_p)_{n_p}(-2|\bx|-|\Bell|_1^{p}+|\bn|_1^{p-1}+|\bx|_1^{p}-\omega)_{n_p}}
{(-\ell_p)_{n_p}(1)_{x_p}(2|\bx|+|\Bell|_1^{p-1}-|\bn|_1^{p-1}-|\bx|_1^{p-1}+\omega+1)_{\ell_p-x_p}}
\\
&\quad\times
\frac{(|\bx|_1^p-|\bn|_1^{p-1}+|\Bell|+a_p+\omega+1)_{\ell_p-x_p}}
    {(-|\bx|_1^{p-1}+|\bn|_1^{p-1}-|\Bell|-\ell_p-a_p-\omega)_{n_p}}
    \frac{(|\bi|_1^p-|\bn|_1^{p-1}+|\Bell|-a_p-\ell_p+\omega^\star)_{\ell_p-i_p}}
    {(-|\bi|_1^{p-1}+|\bn|_1^{p-1}-|\Bell|+a_p+1-\omega^\star)_{n_p}}\\
&\quad\times
\frac{(i_p-\ell_p)_{n_p}(-1)^{\ell_p}(-2|\bi|+|\bn|_1^{p-1}-|\Bell|_1^{p}+|\bi|_1^{p}-\omega^\star)_{n_p}}
     {(1)_{n_p}(2|\bi|-|\bn|_1^{p-1}+|\Bell|_1^{p-1}-|\bi|_1^{p-1}+\omega^\star+1)_{\ell_p-i_p}}\,.
\end{aligned}
\end{align}
Finally, with the help of the shift operator, previous relation simplifies to
\begin{align}\label{eq:cobi_prod4}
\begin{aligned}
U_{\bi}(\bx)&=\!\sum_{\bn=\mathbf{0}}^{\min(\Bell-\bx,\Bell-\bi)}
{\overrightarrow{\prod}}_{p=1}^{N}\Bigg(
\frac{(-\ell_p)_{x_p}(x_p-\ell_p)_{n_p}(-2|\bx|-|\Bell|_1^{p}+|\bx|_1^{p}-\omega)_{n_p}}
{(-\ell_p)_{n_p}(1)_{x_p}(2|\bx|+|\Bell|_1^{p-1}-|\bx|_1^{p-1}+\omega+1)_{\ell_p-x_p}}
\\
&\quad\times
\frac{(|\bx|_1^p+|\Bell|+a_p+\omega+1)_{\ell_p-x_p}}
    {(-|\bx|_1^{p-1}-|\Bell|-\ell_p-a_p-\omega)_{n_p}}
    \frac{(|\bi|_1^p+|\Bell|-a_p-\ell_p+\omega^\star)_{\ell_p-i_p}}
    {(-|\bi|_1^{p-1}-|\Bell|+a_p+1-\omega^\star)_{n_p}}\\
&\quad\times
\frac{(i_p-\ell_p)_{n_p}(-1)^{\ell_p}(-2|\bi|-|\Bell|_1^{p}+|\bi|_1^{p}-\omega^\star)_{n_p}}
     {(1)_{n_p}(2|\bi|+|\Bell|_1^{p-1}-|\bi|_1^{p-1}+\omega^\star+1)_{\ell_p-i_p}}e^{-n_p(\partial_{x_p}+\partial_{i_p})}\Bigg)\cdot\mathbf{1}(x_N)\mathbf{1}(i_N)\,,
\end{aligned}
\end{align}
which proves the result of the theorem.
\endproof

\begin{theo}[\textbf{Biorthogonality}]
The functions $T_{\bi}(\bx)$ and their biorthogonal partner $U_{\bi}(\bx)$
obey the following biorthogonality relation
\begin{align}\label{eq:biortho}
\sum_{\bx=\mathbf{0}}^{\Bell} T_{\bi}(\bx)\,U_{\bj}(\bx)=\delta_{\bi,\bj}\,.
\end{align}
\end{theo}
\begin{exam}
For $N=1$, one gets
\begin{align}\label{eq:cobi_prod4bis}
\begin{aligned}
U_{i}(x)&=(-1)^{\ell}
\frac{(-\ell)_{x}(x+\ell+a+\omega+1)_{\ell-x}
      (i-a+\omega^\star)_{\ell-i}}
     {(1)_{x}(2x+\omega+1)_{\ell-x}(2i+\omega^\star+1)_{\ell-i}}\\
&\quad\times
\pFq{4}{3}{i-\ell,\; x-\ell,\; -x-\ell-\omega,\; -i-\ell-\omega^\star}
          {-\ell,\; -2\ell-a-\omega,\; -\ell+a+1-\omega^\star}{1}\,.
\end{aligned}
\end{align}
Using Whipple's transformation twice, it becomes
\begin{align}
\begin{aligned}
U_i(x)&= (norm)~
\pFq{4}{3}{-i,\; i+\omega^\star,\; -x,\; x+\omega}
          {-\ell,\; -a+\omega^\star,\; \ell+a+\omega+1}{1}\,,\\
(norm)&=
(-1)^x\frac{(-\ell)_x}{(1)_x}
\frac{(1-a+\omega^\star)_{\ell-i}(i-\ell-a-\omega+\omega^\star)_{\ell-i}
      (i+a+1)_{\ell-i}}
     {(2i+\omega^\star+1)_{\ell-i}(-2\ell-a-\omega)_{\ell-i}
      (-\ell+a-\omega^\star+1)_{\ell-i}}\\
&\quad\times
\frac{(x+\ell+a+\omega+1)_{\ell-x}(-x+a-\omega^\star+1)_x
      (-x-\ell-a-\omega)_x}
     {(2x+\omega+1)_{\ell-x}(a+\omega-\omega^\star+1)_x(-\ell-a)_x}
\end{aligned}
\end{align}
and we recognize this as Racah polynomials (up to some normalization).
The biorthogonality relation \eqref{eq:biortho} then becomes the regular
orthogonality relation of the Racah polynomials.
\end{exam}


\section{Limiting cases}
\label{sec:limit}

Tridiagonal pairs of type II are divided in three categories: the first, second and third kind.
The kind of a tridiagonal pair is determined according to the degrees
of the eigenvalues $\theta_i$, $\theta^\star_i$.
The first kind is the most general case, and the second and third kinds can be reached from it
through successive limits.
We detail this below.

\subsection{Racah-type case (First kind)}
This is the general case that was introduced in this paper.
For $i=0,1,\dots,|\Bell|$, the eigenvalues of the first kind can be written as follows:
\begin{align}\label{eq:theta_t1}
 \theta_i=\theta_0 + hi(i+\omega)\,,\qqquad
 \theta^\star_i=\theta^\star_0 + h^\star i(i+\omega^\star)\,,
\end{align}
for some free constants
$\theta_0,\theta_0^\star,h,h^\star,\omega,\omega^\star\in \mathbb{C}$.

The special functions $T_\bi(\bx)$ are given in \eqref{eq:cob_fact}
and resemble a nested product of Racah polynomials.

\subsection{Hahn-type case (Second kind)}
For $i=0,1,\dots,|\Bell|$, the eigenvalues of the second kind can be written as follows:
\begin{align}\label{eq:theta_t2}
 \theta_i=\theta_0 + hi(i+\omega)\,,\qqquad
 \theta^\star_i=\theta^\star_0 + h^\star i\,,
\end{align}
for some free constants
$\theta_0,\theta_0^\star,h,h^\star,\omega\in \mathbb{C}$.
The second kind case is obtained from the first kind case by rewriting
$h^\star\mapsto h^\star t$, $\omega^\star\mapsto1/t$ and letting $t\to0$.

The special functions giving the change of basis in this case become
\begin{align}
\begin{aligned}
 T_{\mathbf{i}}(\bx)=
 \frac{(-1)^{|\mathbf{i}|}}{(|\bx|+\omega)_{|\bx|}}
 {\overrightarrow{\prod}}_{p=1}^{N}
 \hah{i_p,\; x_p,\; |\bx|+\omega}
 {\ell_p,\; |\bx|_{1}^{p-1}+|\Bell|_{p}^{N}+\omega+a_p+1}
 {e^{\partial_{x_p}}}\cdot\mathbf{1}(x_N)\,,
\end{aligned}
\end{align}
with
\begin{align}
\hah{i,\; x,\; a}{\ell,\; b}{Z}
=\binom{\ell}{i}(b)_x\
\pFq{3}{2}{-i,\; -x,\; a}{-\ell,\; b}{Z}\,.
\end{align}
These functions resemble a nested product of Hahn polynomials.

\subsection{Krawtchouk-type case (Third kind)}
For $i=0,1,\dots,|\Bell|$, the eigenvalues of the third kind can be written as follows:
\begin{align}\label{eq:theta_t3}
 \theta_i=\theta_0 + hi\,,\qqquad
 \theta^\star_i=\theta^\star_0 + h^\star i\,,
\end{align}
for some free constants
$\theta_0,\theta_0^\star,h,h^\star\in \mathbb{C}$.
The third kind case is obtained from the second kind case by rewriting
$h\mapsto ht$, $\omega\mapsto1/t$ and letting $t\to0$.

The special functions giving the change of basis in this case become
\begin{align}
\begin{aligned}
 T_{\mathbf{i}}(\bx)=
 \prod_{p=1}^{N}
 \frac{(-\ell_p)_{i_p}}{(1)_{i_p}}
 \pFq{2}{1}{-i_p,\; -x_p}{-\ell_p}{1}
\end{aligned}
\end{align}
This is a simple product of Krawtchouk polynomials
$K_{i_p}(x_p;\phi;N)$ with parameter $\phi=1$.

\section{Outlooks}
\label{sec:outlooks}

In this paper, we provided an explicit expression for the change of basis between
the two bases in the definition of tridiagonal pairs of type II.
It involves special functions which generalize Racah polynomials.
It would be interesting to study if these special functions possess additional
properties such as contiguity relations, duality, etc.

The monovariate Racah polynomials have been used extensively and appear in
numerous contexts.
To cite a few examples,
they provide an explicit expression of the $6j$-symbol introduced in the problem
of the recoupling of three spins in quantum mechanics \cite{Wilson1978},
they characterize $P$- and $Q$- polynomial association scheme
\cite{Leonard1982, BannaiIto1984},
they are closely related to the generic second-order superintegrable model on
the 2-sphere \cite{GenestVinetetal2014},
they define interesting models of free fermions \cite{CrampeNepomechieetal2019}
and also describe exactly solvable birth and death processes \cite{Sasaki2022}.
Our generalization of the Racah polynomials introduced in this paper could also
have similar uses in these different contexts and we are exploring such
potential applications.

In addition, the bispectrality property of the Racah polynomials leads to the
definition of the eponym algebra, the \textit{Racah algebra} (see
\textit{e.g.}~\cite{CrampeGaboriaudetal2021} for a recent review).
This corresponds to the case $N=1$ in this paper.
For $N=1$, the operators $A$ and $A^\star$ satisfy the relations of
the Racah algebra which are more constrained than the ones of the tridiagonal
algebra (see \eqref{eq:TDrelations}).
A detailed understanding of the algebraic properties of the functions for
$N\geq2$ is the next step.
It would be interesting to obtain supplementary relations---in addition to the
ones of the tridiagonal algebra---satisfied by $A$ and $A^\star$ for any given
shape $(\rho_0,\rho_1,\dots,\rho_d)$.
This would provide an algebraic interpretation of the shape of tridiagonal
pairs.

Clearly, results similar to the ones obtained in this paper ought to exist
for tridiagonal pairs of type I or type III and we plan to pursue these in
upcoming work.
As mentioned in the introduction, a change of basis for tridiagonal pairs of
type I has been already computed in \cite{BaseilhacVinetetal2017}. To be
precise, the irreducibility of these pairs was not discussed in this latter
cited paper and cannot be called tridiagonal pairs. However, the constraints on
the parameters such that they become really tridiagonal pairs can be computed
thanks to \cite{Ito2014}.
It would be very interesting to understand the connection between the change of
basis obtained in \cite{BaseilhacVinetetal2017} and the one obtained by
generalizing the results of the present paper. Both are not necessarily
identical since the bases where $A$ or $A^\star$ are diagonal are not unique.
To the best of our knowledge, the computation of the change of basis for
tridiagonal pairs of type III does not appear in the literature.

\subsection*{Acknowledgments}
The authors thank P.~Baseilhac for his interest on this work.
They also thank T.~Ito and P.~Terwilliger for helpful discussions.
N.~Crampé is partially supported by the IRP AAPT of CNRS and
thanks Kyoto University for hospitality.
During the course of this work,
J.~Gaboriaud was an International Research Fellow of
Japan Society for the Promotion of Science
(Postdoctoral Fellowships for Research in Japan (Standard)).
J.~Gaboriaud and S.~Tsujimoto are supported in part by
JSPS KAKENHI Grant Numbers 22F21320 and 22KF0189.
The research of S.~Tsujimoto is also supported by
JSPS Grant-in-Aid for Scientific Research (B), 19H01792 and 24K00528.


\appendix
\section{Notations}\label{sec:notations}
\paragraph{$N$-tuple.}
Boldface letters are used to denote $N$-tuples throughout the paper.
We let $\be_p$ denote the unit $N$-tuples
\begin{align}
 \be_1=(1,0,\dots,0)\,,\qquad
 \be_2=(0,1,0,\dots,0)\,,\qquad
 \dots\,,\qquad
 \be_N=(0,\dots,0,1)\,,
\end{align}
and let $\mathbf{0}$ denote the zero $N$-tuple $\mathbf{0}=(0,\dots,0)$.
The addition of $N$-tuples is also defined by
\begin{align}
 \bn\pm\mathbf{m}=(n_1\pm m_1,\dots,n_N\pm m_N)\,.
\end{align}
We also introduce the notation
\begin{align}
 |\bn|_j^k=\sum_{p=j}^k n_p\,,
 \qqquad \text{with}\quad|\bn|_j^k=0\quad\text{if $j>k$}\,,
\end{align}
as well as the convention
$|\bn|\equiv|\bn|_1^N$.

\paragraph{Hypergeometric functions.}
The generalized hypergeometric function is defined by
\begin{align}
 \pFq{r}{s}{a_1,\; a_2,\;\dots,\; a_r}{b_1,\; \dots,\;b_s}{z}
 =\sum_{k=0}^{\infty}
 \frac{(a_1,\;a_2,\;\dots,\;a_r)_k}{(1)_k(b_1,\;\dots,\;b_s)_k} z^k\,,
\end{align}
with the Pochhammer symbol $(x)_k$ given by
\begin{align}
 (x)_k=(x)(x+1)\dots(x+k-1)\,,\qquad
 (x)_0\equiv1\,,
\end{align}
and
\begin{align}
 (a_1,\;\dots,\;a_n)_k=(a_1)_k\dots(a_n)_k\,.
\end{align}


\printbibliography
\end{document}